\documentclass{article}

\usepackage{amsmath,amssymb,latexsym}
\usepackage{mathrsfs}
\usepackage{verbatim}
\usepackage{color}

\newtheorem{defn}{Definition}[section]

\newtheorem{thm}[defn]{Theorem}
\newtheorem{prop}[defn]{Proposition}

\newtheorem{rem}[defn]{Remark}
\newtheorem{example}[defn]{Example}

\newcommand{\ru}{R(U)}

\newcommand{\uinv}{{U^{-1}}}
\newcommand{\mn}{\mathbb N}
\newcommand{\mr}{\mathbb R}

\newcommand{\nulel}{{\bf 0}}

\newcommand{\seqgr[1]}{\{#1_i\}_{i=1}^\infty }
\newcommand{\sech[1]}{{\cap_{s\in\mn_0}{#1}_s}}

\newcommand{\snorm[1]}{{\|\hspace{-.009in}|{#1}\|\hspace{-.009in}|}}

\def\be{\begin{equation}}
\def\ee{\end{equation}}
\def\newin {\in}

\def\v{\vspace{.1in}}

\def\bp{\noindent{\bf Proof. \ }}
\def\ep{\noindent{$\Box$}}
\def\<{\langle}
\def\>{\rangle}

\begin{document}

\title{
Fr\'echet frames, general definition and expansions\\
}

\author{Stevan Pilipovi\'c and Diana T. Stoeva}

\maketitle

\begin{abstract}  
We define an {\it $(X_1,\Theta, X_2)$-frame} with Banach spaces $X_2\subseteq X_1$,
  $\|\cdot\|_1 \leq \|\cdot\|_2$,  and a $BK$-space $(\Theta, \snorm[\cdot])$. 
  Then  by the use of decreasing sequences of Banach spaces 
  $\{X_s\}_{s=0}^\infty$ and of sequence spaces $\{\Theta_s\}_{s=0}^\infty$, we define a general Fr\' echet frame on the Fr\' echet space
  $X_F=\bigcap_{s=0}^\infty X_s$. 
  We give frame expansions of elements of $X_F$ and 
  its dual $X_F^*$, as well of 
  some of the generating spaces of $X_F$ with convergence in appropriate norms. 
  Moreover, we give necessary and sufficient conditions for a general pre-Fr\' echet frame
  to be a general Fr\' echet frame, as well as for the complementedness of the range of the analysis operator
  $U:X_F\to\Theta_F$. 
\end{abstract}



\section{Introduction} 

\footnotetext[1]{This research was supported by Ministry of Science of Serbia, Project 174024. The work of the second author was also partially supported by DAAD and by the Austrian Science Fund (FWF) START-project FLAME ('Frames and Linear Operators for Acoustical Modeling and Parameter Estimation'; Y 551-N13). }

For given Fr\' echet spaces $X_F=\bigcap_{s=0}^\infty X_s$  and $\Theta_F=\bigcap_{s=0}^\infty \Theta_s$ ($X_s$ and
$\Theta_s$ are Banach space and Banach sequence space with the norms $\|\cdot\|_s$ and  $\snorm[\cdot]_s$, respectively), 
in \cite{ps2}, we have determined
conditions on a sequence $\seqgr[g]$, $g_i\newin X_F^*$, which imply the
existence of $\seqgr[f]$, $f_i\newin X_F$, such that every
$f\newin X_F$ and every $g\newin X_F^*$ can be written as
$f=\sum_{i=1}^{\infty} g_i(f) f_i$ and $g=\sum_{i=1}^{\infty}
g(f_i) g_i$.
These conditions are related to the frame inequalities 
$$
\{g_i(f)\}_{i=1}^\infty\in\Theta_F \ \mbox{ and } \
A_s \|f\|_s\leq \snorm[\{g_i(f)\}_{i=1}^\infty]_s\leq
B_k\|f\|_s, \ f\newin X_F, s\in\mn_0.
$$

In the present paper we are concerned with the series expansions via more general sequences 
in Fr\'echet spaces allowing different norms in the inequalities given above, namely,
$$
\{g_i(f)\}_{i=1}^\infty\in\Theta_F \ \mbox{ and } \
A_k \|f\|_{s_k}\leq \snorm[\{g_i(f)\}_{i=1}^\infty]_{k}\leq
B_k\|f\|_{\widetilde{s}_k}, \ f\newin X_F, k\in\mn_0,
$$
 where $\{s_k\}_{k\in\mn_0}$ and $\{\widetilde{s}_k\}_{k\in\mn_0}$ are increasing  subsequences of $ \mn_0$. 
The statements in \cite[Theorem 5.3]{ps2}(b)(c)  give sufficient conditions for an operator defined on $\Theta_F$ to imply series expansions in some of the generating Banach spaces $X_s$. 
In our main theorems 
 we extend these results and determine conditions on $V:\Theta_F\to X_F$ which are necessary and sufficient for the aim of expansions. 
In particular, in Theorem \ref{diffnew} we prove that the existence of a sequence $\seqgr[f]\newin (X_F)^\mn$ which is a $\Theta_k^*$-Bessel sequence for $X_{s_k}^*$ for every $k\in\mn_0$ and gives series expansions in $ X_{\widetilde{s}_k}$ with convergence in $\|\cdot\|_{s_k}$-norm, is equivalent to the existence of an operator  $V:\Theta_F\rightarrow X_F$ so that $V\{g_i(f)\}_{i=1}^\infty=f, \ \forall f\newin X_F,$ and 
 $\|Vd\|_{s_k}\leq C_k \snorm[d]_k$, $\forall d\in\Theta_F, \,\forall k\in\mn_0$ for some constants $C_k$. 
 While the continuity of $V$ implies series expansions in $X_F$ and the above boundedness properties of $V$ 
 imply series expansions in all the spaces $X_{\widetilde{s}_k}$, $k\in\mn_0$, with convergence in $\|.\|_{s_k}$-norm (see Theorem \ref{diffnew}), 
 in Theorem \ref{th2} we prove that the continuity property of $V$ is enough to imply the existence of subsequence of $\{X_{\widetilde{s}_k}\}_{k=1}^\infty$ with series expansions. 
 In general, it is not easy to show that a general pre-Fr\' echet frame is a general Fr\' echet frame. We devote Section \ref{sec4} (Theorem \ref{equivcompl}) to this problem. Several examples in Sections \ref{egff} and \ref{sec4} illustrate our investigations.

The paper is organized as follows. 
 The notation used in the paper are recalled in Subsection \ref{s2p}. In Section \ref{gff} we give the definition of a general pre-Fr\' echet frame.
 Further, we give some statements needed for the main theorems of Section \ref{egff}.  Section \ref{egff} concerns series expansions via general pre-Fr\' echet frames. We determine sufficient conditions for a general pre-Fr\' echet frame to imply series expansions in a Fr\' echet space and its dual, as well as necessary and sufficient conditions for a general pre-Fr\' echet frame to imply series expansions in spaces generating the Fr\' echet space via a sequence with the Bessel properties.
Moreover, in Section \ref{egff} we define a general Fr\' echet frame and give an example of a general Fr\' echet frame for $X_F$ with respect to $\Theta_F$ which is not a Fr\' echet frame for $X_F$ with respect to $\Theta_F$ (according to the definition in \cite{ps2}). 
In Section \ref{sec4} we give necessary and sufficient conditions for a general pre-Fr\' echet frame to be a general Fr\' echet frame.

Concerning the list of references, one can find more information about papers, related to Banach frame expansions, in the bibliography of \cite{ps2}.

\subsection{Preliminaries}\label{s2p}

Throughout the paper, $(X, \|\cdot \|)$ is a Banach space and
$(X^*, \|\cdot \|_{X^*})$ is its dual, $( \Theta, \snorm[\cdot]) $ is
a Banach sequence space and  $( \Theta^*, \snorm[\cdot]_{\Theta^*}) $ is
the dual of $\Theta$. Recall that $\Theta$ is called a {\it
$BK$-space} if the coordinate functionals on
$\Theta$ are continuous. 
The symbol $e_i$ denotes the $i$-th canonical vector
$\{\delta_{ki}\}_{k=1}^\infty$, $i\newin\mn$.
A $BK$-space $\Theta$ is called a {\it
$\lambda$--$BK$-space} ($\lambda\geq 1$) if it contains all the canonical vectors and
\begin{equation}\label{sol}
\snorm[\{c_i\}_{i=1}^n]_\Theta\leq
\lambda\snorm[\{c_i\}_{i=1}^\infty]_\Theta, \
 n\in\mn, \, \{c_i\}_{i=1}^\infty\in \Theta.
\end{equation}

A $BK$-space for which the canonical vectors form a
Schauder basis, is called a {\it $CB$-space}. If $\Theta$ is a
$CB$-space, then the space $\Theta^\circledast :=\{
\{g(e_i)\}_{i=1}^{\infty} : g\in \Theta^* \}$ with the norm
$\snorm[\{g(e_i)\}_{i=1}^{\infty}]_{\Theta^\circledast}:=\|g\|_{\Theta^*}$
is a $BK$-space, isometrically isomorphic to $\Theta^*$ (see
\cite[p.\,201]{KA}). \label{bkthetastar}
From now on, when $\Theta$ is a
$CB$-space, we will always identify $\Theta^*$ with
$\Theta^\circledast$. In the sequel, linear mappings are
called operators. Recall that an operator 
$\mathcal{P}:X\to Z(\subseteq X)$ is called a {\it projection of $X$ onto $Z$} if $R(P)=Z $ and $\mathcal{P}^2=\mathcal{P}$ (equivalently, if $R(P)=Z $ and $P\vert_Z={\rm Id}_Z$).
 The symbol $\mn_0$ denotes the set $\{0,1,2,3,\ldots\}$.

\section{Definition of a general pre-Fr\' echet frame } \label{gff}

We begin with a generalization of a $\Theta$-frame.

\begin{defn} \label{thetafrnew}
Let $(X_i, \|\cdot \|_i)$, $i=1,2$, be Banach spaces such that
 $X_2 \subseteq X_1$, $\|\cdot\|_1 \leq \|\cdot\|_2$, and let
 $(\Theta, \snorm[\cdot])$ be a $BK$-space. The sequence
 $\seqgr[g]\newin {(X_2^*)}^\mn$ is called an {\it $(X_1,\Theta, X_2)$-frame} with bounds 
 $A, B$ if 
 $0<A\leq B <\infty$
  and for every $f\in X_2$,
\begin{equation}  \label{xd1new}
\{ g_i(f) \}_{i=1}^\infty \in \Theta \ \mbox{ and } \
A\|f\|_1 \leq \snorm[ \{g_i(f)\}_{i=1}^\infty ] \leq B \|f\|_2.
\end{equation}
\end{defn}

\v
Note that when $\Theta$ is a $BK$-space, the validity of the upper inequality in (\ref{xd1new}) for every $f\in X_2$ 
(i.e., $\seqgr[g]$ being a $\Theta$-Bessel sequence for $X_2$) 
implies that $g_i$ must be bounded on $X_2$. 
 When $X_1=X_2=X$, then an $(X_1,\Theta, X_2)$-frame becomes a $\Theta$-frame for $X$.

\v 
We give a generalization of \cite[Theorem 3.3(ii)]{pst}.

\begin{prop} \label{p1new}
Let $(X_i, \|\cdot\|_i)$, $i=1,2$, be Banach spaces such that
 $X_2 \subseteq X_1$ and $\|\cdot\|_1 \leq \|\cdot\|_2$. Let
 $(\Theta, \snorm[\cdot])$ be a $\lambda$--$BK$-space,
 $W$ be a dense subset of $X_2$ 
 and $\seqgr[g]\newin {(X_2^*)}^\mn$.
 If (\ref{xd1new}) holds for all $f\in W$, then
$\seqgr[g]$ is an $(X_1,\Theta, X_2)$-frame with bounds $A, \lambda B$.
\end{prop}
\bp
By \cite[Theorem 3.3(i)]{pst}, it follows that 
\begin{equation*}  
\{g_i(f)\}_{i=1}^\infty\in \Theta \ \mbox{ and } \ \snorm[ \{g_i(f)\}_{i=1}^\infty]_{\Theta} \leq \lambda B \|f\|_2, \ \forall x\in X_2.
\end{equation*}
For the lower inequality, take $f\in
X_2\setminus W$ and a sequence $\{f_n\}\newin W^\mn$ such that $f_n\to f$
when $n\to\infty$ in $\|\cdot\|_2$-norm (and hence, in $\|\cdot\|_1$-norm). 
Since $ \displaystyle \lim_{n\rightarrow \infty }\snorm[\{ g_i(f_n)\}_{i=1}^\infty]
=  \snorm[\{
g_i(f)\}_{i=1}^\infty],$
it follows that $A\|f\|_1\leq \snorm[\{g_i(f)\}_{i=1}^\infty].$
\ep

\v
Now we define a general pre-Fr\'echet frame. 
Let $\{Y_s, | \cdot |_s\}_{s\in\mn_0}$ be a sequence of separable
Banach spaces such that \be \label{fx1} \{\nulel\} \neq
\sech[Y]\subseteq \ldots \subseteq Y_2 \subseteq Y_1 \subseteq Y_0
\ee \be  \label{fx2} |\cdot|_0\leq | \cdot |_1\leq | \cdot |_2\leq
\ldots \ee \be \label{fx3} Y_F :=\sech[Y] \;\; \mbox{is dense in}
\;\; Y_s, \;\; s\in\mn_0. \ee

Then $Y_F$ is a Fr\'echet space with the sequence of norms $ |
\cdot |_s, $ $ s\in\mn_0.$
 We will use such sequences in two cases:

1. $Y_s=X_s$ with norm $\|\cdot\|_s, s\in\mn_0;$ 

2. $Y_s=\Theta_s$ with norm $\snorm[\cdot]_s, s\in\mn_0$.

\begin{defn}\label{fframenew}
Let $X_F$ be a Fr\'echet space determined by the separable Banach spaces $X_s$, $s\in\mn_0$, satisfying (\ref{fx1})-(\ref{fx3}), 
and let $\Theta_F$ be a Fr\'echet space determined by the $BK$-spaces $\Theta_s$, $s\in\mn_0$, satisfying (\ref{fx1})-(\ref{fx3}). A sequence
$\seqgr[g]\newin ({X_F^*})^\mn$ is called a {\it general pre-Fr\' echet frame} (in short, {\it general pre-$F$-frame}) {\it for $X_F$ with respect to $\Theta_F$}  if there exist sequences
$\{\widetilde{s}_k\}_{k\in\mn_0}\subseteq \mn_0, \{s_k\}_{k\in\mn_0}\subseteq \mn_0$
 which increase to $\infty$ with the property
$s_k\leq \widetilde{s}_k$,
$ k\newin\mn_0$, and there exist constants $B_k, A_k>0$, $k\newin\mn_0$, satisfying
\begin{equation}\label{fframestarnew}
\{g_i(f)\}_{i=1}^\infty\in\Theta_F \ \mbox{ and } \
A_k \|f\|_{s_k}\leq \snorm[\{g_i(f)\}_{i=1}^\infty]_{k}\leq
B_k\|f\|_{\widetilde{s}_k}, \ f\newin X_F.
\end{equation}
\end{defn}

The above definition
 reduces to  the definition of a pre-$F$-frame in \cite[Def. 2.3]{ps2} if the sequences $\{s_k\}$ and $\{\widetilde{s}_k\}$ coincide.
We give the definition  of a general $F$-frame after Theorem \ref{diffnew}. We will use the names {\it strict pre-$F$-frame} and {\it strict $F$-frame}
in the cases considered in \cite{ps2} (when $\{s_k\}$ and $\{\widetilde{s}_k\}$ coincide).

\begin{rem} Let $\seqgr[g]$ be a general pre-$F$-frame for $X_F$ with respect to $\Theta_F$ according to Definition \ref{fframenew}. 
One can see that every subsequences $\{X_{p_k}\}_{k=1}^\infty$ of $\{X_s\}_{s=1}^\infty$ and $\{\Theta_{q_k}\}_{k=1}^\infty$ of $\{\Theta_s\}_{s=1}^\infty$ have suitable sub-subsequences so that (\ref{fframestarnew}) holds with the same $\{g_i\}_{i=1}^\infty$ and corresponding sub-subsequences of norms. 
\end{rem}

In the sequel, 
when we consider a general pre-$F$-frame $\seqgr[g]$ for $X_F$ with respect to $\Theta_F$, we always assume that 
$X_F$ is determined by the sequence 
$\{X_s, \|\cdot\|_s\}_{s\in\mn_0}$ of Banach
spaces satisfying (\ref{fx1})-(\ref{fx3}),
$\Theta_F$ is determined by a sequence 
$\{\Theta_s, \snorm[\cdot]_s\}_{s\in\mn_0}$  of
$BK$-spaces satisfying (\ref{fx1})-(\ref{fx3}),
and $\seqgr[g]$ fulfills (\ref{fframestarnew}). 

\begin{rem} \label{gi} Let $\seqgr[g]$ be a general pre-$F$-frame for $X_F$ with respect to $\Theta_F$. 
 For every $i\newin\mn$ and every $k\newin\mn_0$, the functional $g_i$ has a unique continues extension on $X_{\widetilde{s}_k}$ which will be denoted by
$g_i^{\widetilde{s}_k}$. By Proposition \ref{p1new}, for every $k\in\mn_0$, 
the sequence $\{g_i^{\widetilde{s}_k}\}_{i=1}^\infty$ is an ($X_{s_k}, \Theta_{k}, X_{\widetilde{s}_k}$)-frame.
Thus, we can consider operators 
\begin{eqnarray}
U_k : X_{\widetilde{s}_k} \to \Theta_k, & & U_k f= \{g_i^{\widetilde{s}_k}(f)\}_{i=1}^\infty, \ k\in\mn_0,\label{operatorus}\\
U: X_F \to \Theta_F, & & Uf=\{g_i(f)\}_{i=1}^\infty.
\label{operatoru}
\end{eqnarray}
Clearly, they are injective and continuous. 
\end{rem}

\begin{prop}\label{ruclosed}
Let $\seqgr[g]\in (X_F^*)^\mn$ be a general pre-$F$-frame for $X_F$ with respect
to $\Theta_F$. Then the following holds.
\begin{itemize}
\item[{\rm (i)}]
The range $\ru$ of the operator $U$, defined
by (\ref{operatoru}), is closed in $\Theta_F$ 
and the inverse operator 
$U^{-1}: R(U)\to X_F$
is continuous.
\item[{\rm (ii)}] 
The existence of a continuous projection $\mathcal{P}$ from $\Theta_F$ onto $\ru$ (i.e. $\ru$ being complemented in $\Theta_F$) is equivalent to the existence of a continuous operator $V:\Theta_F\rightarrow X_F$ so that $V\{g_i(f)\}_{i=1}^\infty=f$ for all $f\newin X_F$.
\end{itemize}
\end{prop}
\bp 
(i) Let $f_n\in X_F, n\in\mn$, and let $\{Uf_n\}_{n=1}^\infty$ converge to some $b=(b_n)_{n=1}^\infty\in\Theta_F$ in $\Theta_F$ 
when $n\to\infty$.
Fix an arbitrary $k\in\mn_0$.  The
lower inequality in (\ref{fframestarnew}) implies that $\{f_n\}_{n=1}^\infty$ converges in
$X_{s_k}$ when $n\to\infty$ 
 and thus, $\{f_n\}_{n=1}^\infty$ converges in $X_F$ to some element $a\in X_F$. 
Furthermore, the upper inequality in (\ref{fframestarnew}) implies that
 $\{Uf_n\}_{n=1}^\infty$ converges to $U a\in\ru$ in $\Theta_F$. Therefore,
$\ru$ is closed in $\Theta_F$. 
The continuity of $U^{-1}$ is easy to see. 

(ii) 
Let $\mathcal{P}$ be a continues projection from $\Theta_F$ onto $\ru$.
This implies that the operator $V$ defined by $V=\uinv \mathcal{P} : \Theta_F\to X_F$ is also continues.
Clearly, $V$ is an extension of $U^{-1}$.

Conversely, let $V:\Theta_F\to X_F$ be continuous and such that $V\{g_i(f)\}_{i=1}^\infty=f$, $\forall f\in X_F$.
Then the operator $\mathcal{P}:=UV$ is a continuous projection from $\Theta_F$ onto $\ru$.
\ep

\begin{rem} \label{runo} 
Note that the assumption $s_k\to\infty$ is essentially used in Proposition \ref{ruclosed} to prove that $\ru$ is closed in $\Theta_F$.
If $\seqgr[g]\newin {(X_2^*)}^\mn$ is  an  $(X_1,\Theta, X_2)$-frame,
then the range of the operator $\widetilde{U}:X_2\to\Theta, \widetilde{U}f:=\{g_i(f)\}_{i=1}^\infty$, is not necessarily closed in $\Theta$. For example, consider $X_1=\ell^q, \Theta=\ell^2, X_2=\ell^p$ for some $1<p<2<q<\infty$. Let $g_i$ be the $i$-th coordinate functional on $\ell^p$, $i\in\mn$. For every $c=\seqgr[c]\in\ell^p$, 
$\|c\|_{\ell^q}\leq \snorm[\{g_i(c)\}_{i=1}^\infty]_{\ell^2}\leq \|c\|_{\ell^p}$, and thus $\seqgr[g]$ is an $(\ell^q,\ell^2, \ell^p)$-frame. Furthermore, $R(\widetilde{U})$ coincides with $\ell^p$ as sets and thus, $R(\widetilde{U})$ is not closed in $\ell^2$. Note that if $\seqgr[g]\newin {(X_2^*)}^\mn$ is  an  $(X_1,\Theta, X_2)$-frame and $R(\widetilde{U})$ is closed in $\Theta$, then $\seqgr[g]$ must satisfy the lower $\Theta$-frame inequality for $X_2$ and thus $\seqgr[g]$ must be a $\Theta$-frame for $X_2$. 
\end{rem}

\section{Expansions} \label{egff}

In this section we are interested in series expansions via general pre-$F$-frames. 
First note that if $\seqgr[g]$ is a general pre-$F$-frame for $X_F$ with respect to $\Theta_F$ such that every $ f\newin X_F$ can be written as
$f=\sum_{i=1}^{\infty} g_i(f) f_i$ with convergence in $X_F$,  then clearly one can define the operator 
$V:\ru (\subseteq \Theta_F)\rightarrow X_F$ by $V\{g_i(f)\}_{i=1}^\infty=f$ and $V$ must be continuous. Below we continue with sufficient (resp. necessary and sufficient) conditions for the existence of series expansions in $X_F$ and in the generating Banach spaces $X_s$.

\begin{thm} \label{diffnew} Let
$\seqgr[g]$ be a general pre-$F$-frame for $X_F$ with respect to $\Theta_F$.
 \begin{itemize}
 \item[{\rm (a)}] Let $\Theta_{s}$, $s\newin\mn_0$, be  $CB$-spaces, and let there exist a continuous operator $V:\Theta_F\rightarrow X_F$ so that $V\{g_i(f)\}_{i=1}^\infty=f$ for all $f\newin X_F$.
 Then there exists a sequence  $\seqgr[f]\newin (X_F)^\mn$
 such that
\begin{eqnarray}
f&=&\sum_{i=1}^{\infty} g_i(f) f_i,  \   f\newin X_F, \
\mbox{(in $X_F$),} \label{freprnew} \\
g&=&\sum_{i=1}^{\infty} g(f_i) g_i,  \  g\in X_F^*, \
\mbox{(in $X_F^*$).} \label{frepr2new}
\end{eqnarray}
 \item[{\rm (b)}] Let $\Theta_{s}$, $s\newin\mn_0$, be  $CB$-spaces. Then the following three statements are equivalent:
 \begin{itemize}
 \item[ $\mathcal{A}_1:$] There exists an operator  $V:\Theta_F\rightarrow X_F$ so that $V\{g_i(f)\}_{i=1}^\infty=f, \ \forall f\newin X_F,$ and 
 for every $k\in\mn_0$ there is a constant $C_k>0$ satisfying  $\|Vd\|_{s_k}\leq C_k \snorm[d]_k$ for all $d\in\Theta_F$. 

 \item[ $\mathcal{A}_2:$] 
  There exists  $\seqgr[f]\newin (X_F)^\mn$ such that for every $k\in\mn_0$, 
  $\seqgr[f]$ is a $\Theta_k^*$-Bessel sequence for $X_{s_k}^*$, and (\ref{freprnew}) holds.

  \item[ $\mathcal{A}_3:$]
 There exists  $\seqgr[f]\newin (X_F)^\mn$ such that for every $k\in\mn_0$, $\seqgr[f]$
 is a $\Theta_k^*$-Bessel sequence for $X_{s_k}^*$,  and
 \begin{eqnarray}  \label{fsreprnew}
f&=&\sum_{i=1}^{\infty} g_i^{\widetilde{s}_k}(f) f_i \ \mbox{in $\|.\|_{s_k}$-norm}, \  f\newin X_{\widetilde{s}_k}.
\end{eqnarray}

\end{itemize}

In particular, 
the equivalent conditions $\mathcal{A}_1$-$\mathcal{A}_3$ imply validity of (\ref{freprnew})-(\ref{fsreprnew}) with a same sequence $\seqgr[f]=\{Ve_i\}_{i=1}^\infty$.

\item[{\rm(c)}] Let $\Theta_s$ and  $\Theta^*_s$, $s\newin\mn_0$, be $CB$-spaces. Then $\mathcal{A}_1$ is equivalent to 
\begin{itemize}
  \item[ $\mathcal{A}_4:$]
 There exists $\seqgr[f]\newin (X_F)^\mn$, such that for every $k\in\mn_0$,
  $\seqgr[f]$ is an $(X^*_{\widetilde{s}_k},\Theta_k^*, X^*_{s_k})$-frame, 
and
\begin{equation} \label{xsnew}
g\vert_{_{X_{\widetilde{s}_k}}}=\sum_{i=1}^{\infty} g(f_i) g_i^{\widetilde{s}_k} 
 \mbox{ in  $\|.\|_{X^*_{\widetilde{s}_k}}$-norm }, \  g\in X_{s_k}^*. 
\end{equation}

\end{itemize}

In particular, 
the equivalent conditions $\mathcal{A}_1$-$\mathcal{A}_4$ imply  validity of (\ref{freprnew})-(\ref{xsnew}) with a same sequence $\seqgr[f]=\{Ve_i\}_{i=1}^\infty$.

\end{itemize}
\end{thm}

\bp 
 (a) Let $f_i=Ve_i, i\in\mn$. 
 By the assumptions, one can write, in $\Theta_F$, $Uf=\sum_{i=1}^\infty g_i(f) e_i$, $f\in X_F$. 
 For every $f\in X_F$, the continuity of $V$ implies that 
$V\left(\sum_{i=1}^n g_i(f)e_i\right)\to VUf=f \ \mbox{in}\ X_F\ \mbox{as}\ n\to\infty$, and this gives  (\ref{freprnew}) and (\ref{frepr2new}).

\vspace{.05in}
(b) $\mathcal{A}_1$ $\Rightarrow$ $\mathcal{A}_2, \mathcal{A}_3$:  
Assume that $\mathcal{A}_1$ holds. Let  $f_i=Ve_i$, $i\in\mn$.
First observe that $V$ is continuous and hence, by (a), the representations 
(\ref{freprnew}) and (\ref{frepr2new}) hold.
Fix an arbitrary $k\in\mn_0$. 
The operator $V$ has a bounded linear extension $V_k :\Theta_k \to X_{s_k}$ and we can
consider the bounded operator $V_k U_k : X_{\widetilde{s}_k} \to X_{s_k}$.
Since $V_k U_k f=f$ for every $f\in X_F$,
$X_F$ is dense in $X_{\widetilde{s}_k}$ and $\|.\|_{s_k} \leq \|.\|_{\widetilde{s}_k}$, it follows that
$V_k U_k f=f$ for every $f\in X_{\widetilde{s}_k}$.
Let $f\in X_{\widetilde{s}_k}$. 
Then $\{g_i^{\widetilde{s}_k}(f)\}_{i=1}^\infty\in \Theta_k$ (see Remark \ref{gi}) 
and
$$ \sum_{i=1}^n g_i^{\widetilde{s}_k} (f) f_i = V_k\left(\sum_{i=1}^n g_i^{\widetilde{s}_k} (f) e_i\right)  
{\underset{n\to\infty}{\xrightarrow{\hspace*{.5cm}}}} \
V_k U_k f=f  \  \mbox{in $\|.\|_{s_k}$}.$$
It is clear that $\seqgr[f]$ is a $\Theta_k^*$-Bessel sequence for $X_{s_k}^*$.

$\mathcal{A}_2$ $\Rightarrow$ $\mathcal{A}_1$:
Assume that $\mathcal{A}_2$ holds. Fix $k\in\mn_0$. Since $\seqgr[f]$ is a $\Theta_k^*$-Bessel sequence for $X_{s_k}^*$, it follows that the synthesis operator $T_k$ given by $T_k\seqgr[d]= \sum_{i=1}^\infty d_i f_i$ is well defined (and bounded) from $\Theta_k$ into $X_{s_k}$ \cite{CCS}. For $\seqgr[d]\in\Theta_F$, the series $\sum_{i=1}^\infty d_i f_i$ converges in $X_{s_k}$ for every $k\in\mn_0$, and thus, it converges in $X_F$. Then we can consider the operator $V:\Theta_F\to X_F$ defined by $V\seqgr[d]= \sum_{i=1}^\infty d_i f_i$. For every $d\in\Theta_F$, 
$$\|Vd\|_{s_k}=\|T_kd\|_{s_k} \leq \|T_k\|\cdot \snorm[d]_k, \ \forall k\in\mn_0.$$
Further, the validity of (\ref{freprnew}) implies that  
$V\{g_i(f)\}_{i=1}^\infty=f$ for all $f\newin X_F$.
Clearly,  $Ve_i=f_i$, $i\in\mn$.

$\mathcal{A}_3$ $\Rightarrow$ $\mathcal{A}_2$: Assume that $\mathcal{A}_3$ holds.
The representations in (\ref{fsreprnew}) imply that 
$f=\sum_{i=1}^{\infty} g_i(f) f_i$  in $\|.\|_{s_k}$-norm for every $k\in\mn_0$ and every $f\newin X_F$, which implies that (\ref{freprnew}) holds. 

\vspace{.05in}
(c) $\mathcal{A}_1$ $\Rightarrow$ $\mathcal{A}_4$:
 Assume that $\mathcal{A}_1$ holds. Let $\seqgr[f]$ be given as in (b) and fix $k\in\mn_0$. 
 Then $\seqgr[f]$ is a $\Theta_k^*$-Bessel sequence for $X_{s_k}^*$. 
   Therefore the synthesis operator  $\widetilde{T}_k$ given by $\widetilde{T}_k \{d_i\}_{i=1}^\infty = \sum_{i=1}^\infty d_i g_i^{\widetilde{s}_k}$ is well defined and bounded from $ \Theta_k^*$ into $X_{\widetilde{s}_k}^*$  (see \cite{CCS}). 
   Let $g\in X_{s_k}^*$. 
   For every $f\in X_{\widetilde{s}_k}$, 
  it follows by (b) that
$\sum_{i=1}^n g_i^{\widetilde{s}_k}(f)f_i\to f$ in $\|\cdot\|_{s_k}$-norm when $n\to\infty$,
which implies that  $g(\sum_{i=1}^n g_i^{\widetilde{s}_k}(f)f_i)\to g(f)$ when $n\to\infty$. 
 Furthermore,
\begin{eqnarray*}
\|g\|_{X_{\widetilde{s}_k}^*} 
&=& \sup_{f\in X_{\widetilde{s}_k}, \|f\|_{\widetilde{s}_k}\leq 1} \left|\sum_{i=1}^\infty g_i^{\widetilde{s}_k}(f) g(f_i)\right| = \sup_{f\in X_{\widetilde{s}_k}, \|f\|_{\widetilde{s}_k}\leq 1} |\widetilde{T}_k\{g(f_i)\}_{i=1}^\infty (f)| \\
&\leq&  
 \|\widetilde{T}_k\|\, \snorm[\{ g(f_i)\}_{i=1}^\infty]_{\Theta_k^*}.
  \end{eqnarray*}
    Therefore, $\seqgr[f]\newin (X_F)^\mn$ is an
 $(X^*_{\widetilde{s}_k},\Theta_k^*, X^*_{s_k})$-frame.
   
   To prove (\ref{xsnew}), denote the canonical basis of $\Theta_k^*$ by $\{\delta_i\}_{i=1}^\infty$.
Let $g\in X_{s_k}^*$. Then $g\vert_{X_{\widetilde{s}_k}} \in X_{\widetilde{s}_k}^*$ and

     \begin{eqnarray*}
 \| g\vert_{_{X_{\widetilde{s}_k}}} - \sum_{i=1}^n g(f_i) g_i^{\widetilde{s}_k}\|_{X_{\widetilde{s}_k}^*} 
   &= & \sup_{f\in X_{\widetilde{s}_k}, \|f\|_{\widetilde{s}_k}\leq 1} 
   | \sum_{i=1}^\infty g_i^{\widetilde{s}_k}(f) g(f_i)  - \sum_{i=1}^n g(f_i) g_i^{\widetilde{s}_k}(f)| \\ 
      &\leq&   
      \|\widetilde{T}_k\| \, |||\sum_{i=n+1}^\infty g(f_i) \delta_i|||_{\Theta_k^*}
      {\underset{n\to\infty}{\xrightarrow{\hspace*{.5cm}}}} \ 0. 
        \end{eqnarray*}

  $\mathcal{A}_4\Rightarrow \mathcal{A}_3$: 
  Assume that $\mathcal{A}_4$ holds. For $k\in\mn_0$, let $B_k$ denote a $\Theta_k^*$-Bessel bound for $\seqgr[f]$. Fix an arbitrary $k\in\mn_0$. For every $f\in X_{\widetilde{s}_k}$, 
    \begin{eqnarray*}
|| f - \sum_{i=1}^n g_i^{\widetilde{s}_k}(f) f_i ||_{s_k} & = & \sup_{g\in
X_{s_k}^*, ||g||_{X_{s_k}^*}=1}
| g(f)- \sum_{i=1}^n g(f_i)g_i(f)| \\
& = & \sup_{g\in X_{s_k}^*, ||g||_{X_{s_k}^*}=1}
|\sum_{i=1}^\infty  g(f_i)g_i(f) - \sum_{i=1}^n
g(f_i)g_i(f)| \\
& = & \sup_{g\in X_{s_k}^*, ||g||_{X_{s_k}^*}=1}
|\sum_{i=n+1}^\infty  g(f_i)g_i(f) | \\
& \le &
\sup_{g\in X_{s_k}^*, ||g||_{X_{s_k}^*}=1}
\snorm[ \{ g(f_i) \}_{i=1}^\infty]_{\Theta_k^*} \ |||\sum_{i=n+1}^\infty g_i(f) e_i |||_{k} \\
& \le & B_k \ |||\sum_{i=n+1}^\infty g_i(f) e_i |||_{k}
 \to  0 \ \mbox{as} \ n \to \infty.
\end{eqnarray*}
Therefore, (\ref{fsreprnew}) holds. 
\ep

\v
Theorem \ref{diffnew} extends \cite[Theorem 5.3]{ps2} 
and improves the formulation of  \cite[Theorem 5.3]{ps2}(a), where it is silently 
assumed that $\Theta_s$, $s\in\mn_0$, are $CB$-spaces.

\v As one can see in Theorem \ref{diffnew}, 
the continuity property of the operator $V$ implies series expansions in $X_F$, while some boundedness properties of $V$ 
 imply series expansions in all the spaces $X_{\widetilde{s}_k}$, $k\in\mn_0$, with convergence in $\|.\|_{s_k}$-norm. 
  Below we prove that the continuity property of $V$ is enough to imply the existence of a subsequence 
 $\{X_{\widetilde{w}_{j}}\}_{j=0}^\infty$ of $\{X_{\widetilde{s}_k}\}_{k=0}^\infty$ such that one has series expansions in $X_{\widetilde{w}_{j}}$, $j\in\mn_0$, with convergence in appropriate norms. 
 
\begin{thm} \label{th2}
 Let $\seqgr[g]$ be a general pre-$F$-frame for $X_F$ with respect to $\Theta_F$ and let $\Theta_{s}$, $s\newin\mn_0$, be  $CB$-spaces. Assume that there exists a continuous operator $V:\Theta_F\rightarrow X_F$ so that $V\{g_i(f)\}_{i=1}^\infty=f$ for all $f\newin X_F$. Then there exist  sequences $\{w_j\}_{j\in\mn_0}$, $\{r_j\}_{j\in\mn_0}$, $\{\widetilde{w}_j\}_{j\in\mn_0}$ which increase to $\infty$ and there exist constants $\widetilde{A}_j, \widetilde{B}_j$, $j\in\mn_0$, such that
 for every $j\in\mn_0$, 
\begin{equation*}
\widetilde{A}_j \|f\|_{w_j}\leq \snorm[\{g_i(f)\}_{i=1}^\infty]_{r_j}\leq
\widetilde{B}_j\|f\|_{\widetilde{w}_j}, \ \forall f\newin X_F.
\end{equation*} 
  Moreover, there exists a sequence  $\seqgr[f]\newin (X_F)^\mn$ such that for every $j\in\mn_0$, $\seqgr[f]$
 is a $\Theta_{r_j}^*$-Bessel sequence for $X_{w_j}^*$ and
 \begin{eqnarray*}  
f&=&\sum_{i=1}^{\infty} g_i^{\widetilde{w}_j}(f) f_i \ \mbox{in $\|.\|_{w_j}$-norm}, \  f\newin X_{\widetilde{w}_{j}}.
\end{eqnarray*}

\end{thm}

\bp 
Assume that $V:\Theta_F\rightarrow X_F$ is a continuous operator satisfying $V\{g_i(f)\}_{i=1}^\infty=f, \ \forall f\newin X_F$. Then  
\begin{equation} \label{vb}
\forall k\in\mn_0, \   \exists  p_k\in\mn_0  \mbox{ and } \ \exists C_{k}  
\mbox{ so that } \|Vd\|_{s_k}\leq C_{k} \,\snorm[d]_{p_k}, \forall d\in \Theta_F.
\end{equation} 
Consider the sequence $\{n_k\}_{k\in\mn_0}$ defined by
$n_k=\max (k, p_k)$, $k\in\mn_0$. 
Clearly, the sequence  $\{n_k\}_{k\in\mn_0}$ is not bounded. 
Take a strictly increasing subsequence $\{n_{k_j} \}_{j\in\mn_0}$ of $\{n_k\}_{k\in\mn_0}$. 
Let $q\in\{k_j\}_{j=1}^\infty$. Since $n_q\geq q$, (\ref{fframestarnew}) implies 
\begin{equation*}
A_q \|f\|_{s_q}
\leq  
\snorm[\{g_i(f)\}_{i=1}^\infty]_{n_q}
\leq
B_{n_q}\|f\|_{\widetilde{s}_{n_q}}, \ \forall f\newin X_F.
\end{equation*}
Since $n_q\geq p_q$,  (\ref{vb}) implies that  $\|Vd\|_{s_q}\leq C_q \,\snorm[d]_{n_q}$, $\forall d\in \Theta_F$.
 Now Theorem \ref{diffnew}(b) implies that there exists  $\seqgr[f]\newin (X_F)^\mn$ such that for every $q\in\{k_j\}_{j=1}^\infty$, $\seqgr[f]$
 is a $\Theta_{n_q}^*$-Bessel sequence for $X_{s_q}^*$ and
 \begin{eqnarray*}  
f&=&\sum_{i=1}^{\infty} g_i^{\widetilde{s}_{n_q}}(f) f_i \ \mbox{in $\|.\|_{s_q}$-norm}, \  f\newin X_{\widetilde{s}_{n_q}}.
\end{eqnarray*}

For $j\in\mn_0$, take $\widetilde{A}_j=A_{k_j}$, $\widetilde{B}_j=B_{n_{k_j}}$, 
 $w_j=s_{k_j}$, $r_j=n_{k_j}$, $\widetilde{w}_j=\widetilde{s}_{n_{k_j}}$. 
\ep

\v
Motivated by 
Theorems \ref{diffnew} and \ref{th2},
we give the definition of a general $F$-frame:

\begin{defn}\label{fframe}
 The sequence $\seqgr[g]\in (X_F^*)^\mn$ is called a {\it general Fr\' echet frame} (in short, {\it general $F$-frame}) {\it for $X_F$ with respect to $\Theta_F$} if $\seqgr[g]$ is a general pre-$F$-frame for $X_F$ with respect to $\Theta_F$ and 
there exists a continuous operator $V:\Theta_F\rightarrow X_F$ so that $V\{g_i(f)\}_{i=1}^\infty=f$ for all $f\newin X_F$. 
\end{defn}

We end the section with an example of a general $F$-frame for $X_F$ with respect to $\Theta_F$ which is not a strict $F$-frame for $X_F$ with respect to $\Theta_F$.

\begin{example} \label{exf1}
{\rm
Let $A$ be a self-adjoint differential operator (for example, one dimensional  normalized harmonic oscillator $(-d^2/dx^2+1)/2$) 
with eigenvalues $\lambda_j=j, j\in\mathbb N,$ and eigenfunctions $\psi_j, j\in\mathbb N$ (Hermite functions)
which make an orthonormal basis of $X_0=L^2(\mathbb R).$ 
For $s\in\mn$, let $X_s$ be the Hilbert space consisting of $L^2-$functions $\phi=\sum_{j=1}^{\infty}a_j\psi_j, a_j\in\mathbb C, j\in\mathbb N,$ with the property
 $\sum_{j=1}^\infty |a_j|^2 j^{2s}<\infty$ and  with the inner product
$$\langle\phi_1,\phi_2\rangle_s=\sum_{j=1}^{\infty}a_{1,j}\overline{a_{2,j}}j^{2s}.
$$
Then $X_F$ is the Fr\'echet space $\mathcal{S}(\mr)$, the Schwartz class of rapidly decreasing functions, 
and its dual is $X_F^*=\mathcal{S}'(\mr)$, the space of tempered distributions. 
For the sequence spaces $\Theta_s$, $s\in\mn_0,$ we take
$$\{d_j\}_{j=1}^\infty\in\Theta_s \mbox{ if and only if } \sum_{j=1}^\infty|d_j|^2 j^{2s}<\infty,
$$
with the usual inner product; $\Theta_F$ is the space of rapidly decreasing sequences. 
Note, 
the space $\Theta_F$ defined above is actually the 
space of the type
\begin{equation*}
\Lambda_\infty(\alpha) = \left\{ \{x_j\}_{j=1}^\infty \ \ : \ \ \snorm[\{x_j\}_{j=1}^\infty]_s:=
\left(\sum_{j=1}^\infty |x_j|^2 e^{2s\alpha_j} \right)^{1/2}<\infty, \ \forall s\in\mn\right\},
\end{equation*} 
with $\alpha_j=\log j$,  
$j\in\mn$. For more information about the spaces $\Lambda_\infty(\alpha)$ we refer to  \cite[Sect. 29]{MV}.

Let $r\in\mn$ be given and let $\{b_j\}_{j=1}^\infty$ be a
sequence of complex numbers such that  
\begin{equation*}
 |b_j| = \left\{  \begin{array}{rl}  1, & j=1, 3, 5,   \ldots;\\
j^r, & j=2, 4, 6, \ldots. 
\end{array}
\right.
\end{equation*}
Let $g_j=b_j\psi_j, j\in\mathbb N,$ and $\phi=\sum_{j=1}^\infty a_j \psi_j\in X_F$. 
Then $\{g_j(\phi)\}_{j=1}^\infty\in\Theta_F$ and
$$||\phi||_s\leq |||\{g_j(\phi)\}_{j=1}^\infty|||_s =\sqrt{ \sum_{j=1}^\infty |a_jb_j|^2 j^{2s}}\leq ||\phi||_{s+r}, s\in \mn_0.
 $$
 Thus,    $\{g_j\}_{j=1}^\infty$  is a general pre-$F$-frame for $X_F$ with 
 respect to $\Theta_F$. 
 Let the operator $U$ be given by (\ref{operatoru}). Observe that $R(U)=\Theta_F$. 
  Therefore,
$\{g_j\}_{j=1}^\infty$  is a general $F$-frame for $X_F$ with respect to $\Theta_F$.

\v
 Furthermore, 
  we will show that $\{g_j\}_{j=1}^\infty$  is not a strict pre-$F$-frame for $X_F$ with respect to $\Theta_F$. 
 Conversely, 
  assume that there exist constants $A_s\in(0,\infty)$, $B_s\in(0,\infty)$, $s\in\mn_0$, and a sequence $\{n_s\}_{s=0}^\infty$,  
satisfying  
\begin{equation}\label{str}
A_s||\phi||_{n_s}\leq  |||\{g_j(\phi)\}_{j=1}^\infty|||_s \leq B_s ||\phi||_{n_s}, \forall \phi\in X_F.
\end{equation}
Fix an arbitrary $s\in\mn_0$. 

If $n_s\leq s$, then (\ref{str}) applied to $\psi_j$, $j\in\mn$, implies that 
$|b_j|\leq B_s j^{n_s-s} \leq B_s$ for all $j\in\mn$, which leads to a contradiction.

If $n_s>s$, then (\ref{str}) applied to $\psi_j$, $j$ - odd, implies that 
$ A_s  \leq j^{s-n_s}$ for all odd $ j$,
which leads to a contradiction.

Therefore, (\ref{str}) can not hold.

\v
Note that if the sequence $\{b_j\}_{j=1}^\infty$ is defined by   
$ |b_j| =  j^r, j\in\mn,$
then $|||\{g_j(\phi)\}_{j=1}^\infty|||_s = ||\phi||_{s+r}$, $s\in \mn_0$, $\phi\in X_F$,
and the sequence $\{g_j\}_{j=1}^\infty$  is a strict pre-$F$-frame for $X_F$ with respect to $\Theta_F$.
}
\end{example} 

\section{On the existence of a continuous projection from $\Theta_F$ onto $\ru$} \label{sec4}

Let $\{g_i\}_{i=1}^\infty$ be a general pre-$F$-frame for $X_F$ with respect to $\Theta_F$. 
As it is shown in Section \ref{egff}, 
the existence of a continuous operator $V:\Theta_F\rightarrow X_F$ such that $V\{g_i(f)\}_{i=1}^\infty=f$ for all $f\newin X_F$ is important for series expansions in $X_F$ 
and 
 in some of the generating Banach spaces 
(see Theorems \ref{diffnew} and \ref{th2}).
Here we consider equivalences of this condition. 
 Thus, we give necessary and sufficient conditions for a general pre-Fr\' echet frame
  to be a general Fr\' echet frame.

\begin{thm} \label{equivcompl} 
Let $\seqgr[g]$ be a general pre-$F$-frame for $X_F$ with respect to $\Theta_F$
and let $\Theta_{s}$, $s\newin\mn_0$, be  $CB$-spaces. Then the following statements are equivalent.
\begin{itemize}

\item[{\rm (i)}] There exists a continuous operator $V:\Theta_F\rightarrow X_F$ so that $V\{g_i(f)\}_{i=1}^\infty=f$ for all $f\newin X_F$. 
\item[{\rm (ii)}]  There exists  $\seqgr[f]\newin (X_F)^\mn$ such that $\sum_{i=1}^\infty c_i f_i$ converges in $X_F$ for every $\{c_i\}_{i=1}^\infty \in\Theta_F$ and (\ref{freprnew}) holds.
\item[{\rm (iii)}] There exist  $\seqgr[f]\newin (X_F)^\mn$ and sequences $\{w_j\}_{j\in\mn_0}$, $\{r_j\}_{j\in\mn_0}$, such that 
  $\seqgr[f]$  is a $\Theta_{r_j}^*$-Bessel sequence for $X_{w_j}^*$ for every $j\in\mn_0$ and (\ref{freprnew}) holds. 
   \end{itemize}
\end{thm}
\bp 
The proof is similar to the one of \cite[Prop. 3.4]{CCS}, extending it to the Fr\'echet case.

\sloppy   (i) $\Rightarrow$ (ii): For $i\in\mn$, define $f_i=Ve_i$. For every $\{c_i\}_{i=1}^\infty\in\Theta_F$, 
   $V(\sum_{i=1}^n c_i e_i)\to V(\{c_i\}_{i=1}^\infty)$ in $X_F$ as $n\to\infty$, which implies that $\sum_{i=1}^\infty c_i f_i$ is convergent in $X_F$. Furthermore, for every $f\in X_F$,
   $ f=V\{g_i(f)\}_{i=1}^\infty=\sum_{i=1}^\infty g_i(f) f_i.$

   (ii) $\Rightarrow$ (i): Assume that (ii) holds and consider the operator $V:\Theta_F\to X_F$ defined by $V(\{c_i\}_{i=1}^\infty)=\sum_{i=1}^\infty c_i f_i$, $\{c_i\}_{i=1}^\infty\in\Theta_F$. 
   Fix an arbitrary $N\in\mn$ and consider the operator $V_N:\Theta_F\to X_F$ defined by $V_N(\{c_i\}_{i=1}^\infty)=\sum_{i=1}^N c_i f_i$, $\{c_i\}_{i=1}^\infty\in\Theta_F$. For every $k\in\mn$, denote the $i$-th coordinate functional on $\Theta_k$ by $E^k_i$ and observe that
   for every $\{c_i\}_{i=1}^\infty\in\Theta_F$ one has
     \begin{eqnarray*}
   \|V_N(\{c_i\}_{i=1}^\infty)\|_{s_k} &=& \|\sum_{i=1}^N c_i f_i\|_{s_k}\leq \sum_{i=1}^N |c_i|\cdot \| f_i\|_{s_k}  \\ 
    &\leq&  \sum_{i=1}^N \| E^k_i\|\cdot \snorm[\{c_i\}_{i=1}^\infty]_k \cdot \| f_i\|_{s_k}\\
    &=& \left(\sum_{i=1}^N \| E^k_i\|\cdot  \| f_i\|_{s_k}\right) \snorm[\{c_i\}_{i=1}^\infty]_k,
   \end{eqnarray*} 
which implies that $V_N$ is continuous on $\Theta_F$. Now the Principle of Uniform Boundedness (see \cite[II.1.17]{DS}) implies that $V$ is continuous. Furthermore, for every $f\in X_F$, $V\{g_i(f)\}_{i=1}^\infty=\sum_{i=1}^\infty g_i(f) f_i=f$.

(i) $\Rightarrow$ (iii): By the proofs of Theorems \ref{th2} and \ref{diffnew}(a)(b), it follows that the sequence $f_i=Ve_i$, $i\in\mn$, fulfills the required properties. 

(iii) $\Rightarrow$ (ii):  Let $\seqgr[d]\in\Theta_F$. 
Similar to the proof of ($\mathcal{A}_2$ $\Rightarrow$ $\mathcal{A}_1$) in Theorem \ref{diffnew}, 
 the series $\sum_{i=1}^\infty d_i f_i$ converges in $X_{\omega_j}$ for every $j\in\mn_0$, and thus, it converges in $X_F$. 
\ep

\v 
Recall that Proposition \ref{ruclosed}(ii) contains one more equivalent condition of Theorem \ref{equivcompl}(i), namely, the existence of a continuous projection of $\Theta_F$ onto $\ru$. 
In Example \ref{exf1} we constructed a general pre-$F$-frame with $\ru=\Theta_F$ and thus it was automatically a general $F$-frame. 
In  Example \ref{exf2} below we construct a general pre-$F$-frame with $\ru\subsetneq \Theta_F$ and show 
the existence of a continuous projection of $\Theta_F$ onto $\ru$.

\begin{example} \label{exf2}
{\rm 
Let $\psi_i$, $i\in\mn$, and $\Theta_s$, $s\in\mn_0$, be defined as in Example \ref{exf1}. 
Let $X_0$ be the closed linear span of the functions $\psi_{2k}$, $k\in\mn$, in $L^2(\mr)$. 
For $s\in\mn$, define $X_s$ to be the space
of those $L^2-$functions $\phi=\sum_{j=1}^{\infty}a_j\psi_{2j}$ ($\in X_0$), $a_j\in\mathbb C, j\in\mathbb N,$ with the property
 $\sum_{j=1}^\infty |a_j|^2 (2j)^{2s}<\infty$ and  with the inner product
$$\langle\phi_1,\phi_2\rangle_s=\sum_{j=1}^{\infty}a_{1,j}\overline{a_{2,j}}(2j)^{2s}.
$$  
Then $\{X_s\}_{s\in\mn_0}$  is a sequence of Hilbert spaces which satisfies (\ref{fx1})-(\ref{fx3}). 

Let $r\in\mn$  be given and let the sequence $\{b_j\}_{j=1}^\infty$ be defined by
$$b_1=b_2=1, b_3=b_4=4^r, b_5=b_6=1, b_7=b_8=8^r, \ldots\,.$$ 
Define
$$\{g_j\}:=\{b_1 \psi_2, b_2 \psi_2, b_3 \psi_4, b_4 \psi_4, b_5 \psi_6, b_6 \psi_6, \ldots \}.$$
Let $\phi=\sum_{j=1}^\infty a_j \psi_j\in X_F$. 
Then $\{g_j(\phi)\}_{j=1}^\infty\in\Theta_F$ and
$$||\phi||_s\leq |||\{g_j(\phi)\}_{j=1}^\infty|||_s \leq\sqrt{ 2\sum_{j=1}^\infty |a_jb_{2j}|^2 (2j)^{2s}}\leq \sqrt{2}||\phi||_{s+r}, s\in \mathbb N_0.
 $$
 Thus,    $\{g_j\}_{j=1}^\infty$  is a general pre-$F$-frame for $X_F$ with 
 respect to $\Theta_F$.

Let the operator $U$ be given by (\ref{operatoru}).
Clearly, 
 $R(U)\subsetneq \Theta_F$. 
We will prove that $\ru$ is complemented in $\Theta_F$. 
Consider the operator $P$ defined on $\Theta_F$ by
$$P(\{d_j\}_{j=1}^\infty):=\{d_2, d_2, d_4, d_4, d_6, d_6, \ldots \}, \ \{d_j\}_{j=1}^\infty \in\Theta_F.$$
Fix an arbitrary $\{d_j\}_{j=1}^\infty \in\Theta_F$. Let 
$$\{a_j\}_{j=1}^\infty:=\left\{d_2, \frac{d_4}{4^r}, d_6, \frac{d_8}{8^r}, d_{10}, \frac{d_{12}}{12^r},\ldots \right\} \ \mbox{ and } \ \phi=\sum_{j=1}^\infty a_j \psi_{2j} \ \mbox{in $L^2$}.$$ 
Clearly, $\phi\in X_F$ and $P\{d_j\}_{j=1}^\infty =\{g_j(\phi)\}_{j=1}^\infty$. 
 Hence, $R(P)\subseteq \ru$. 
For every $\phi\in X_F$, one has $ U\phi\in\Theta_F$ and $U\phi=P(U\phi)$. 
Therefore, $\ru=R(P)$ and $P(=P^2)$ is a projection of $\Theta_F$ onto $\ru$. 
Further on, for every $\{d_j\}_{j=1}^\infty \in\Theta_F$,
$$ \snorm[P(\{d_j\}_{j=1}^\infty)]_s^2\leq 2 \sum_{j=1}^\infty |d_{2j}|^2 (2j)^{2s}\leq 2 \snorm[\{d_j\}_{j=1}^\infty]_2.$$
Thus, $P$ is continuous and 
$\{g_j\}_{j=1}^\infty$  is a general $F$-frame for $X_F$ with respect to $\Theta_F$.

 Furthermore, 
 in a similar way as in Example \ref{exf1} one can 
 show that $\{g_j\}_{j=1}^\infty$  is not a strict pre-$F$-frame for $X_F$ with respect to $\Theta_F$. 

}
\end{example}

\v In Example \ref{exf2}, the space $\Theta_F$ is of the type $\Lambda_\infty(\alpha)$ and $\ru$ is closed in $\Theta_F$ by Proposition \ref{ruclosed}. 
In this example we prove that $\ru$ is complemented in $\Theta_F$ by a direct construction of a continuous projection of $\Theta_F$ onto $\ru$. 
For another necessary and sufficient condition for $\ru$ to be complemented in $\Lambda_\infty(\alpha)$, we refer to \cite[Section 30, Exercise 2]{MV}:

{\it
\lq\lq A closed subspace $E$ of $\Lambda_\infty(\alpha)$ is complemented in $\Lambda_\infty(\alpha)$ if and only if $E$ has the property $(\Omega)$ and $\Lambda_\infty(\alpha)/E$ has the property (DN).\rq\rq}

\noindent For the definitions of the properties $(\Omega)$  and (DN), see \cite[Sect. 29]{MV}.

\v
 {\bf Acknowledgments} The second author is grateful for the hospitality of the University of Novi Sad, where most of the work on the paper was done.

\v\noindent
S. Pilipovi\'c\\
Department of Mathematics and Informatics\\
University of Novi Sad\\
Trg D. Obradovi\'ca 4\\
21000 Novi Sad, Serbia\\
stevan.pilipovic@dmi.uns.ac.rs

\v\noindent
D.\,T. Stoeva\\
Acoustics Research Institute \\
Wohllebengasse 12-14\\
 Vienna A-1040, Austria\\
  dstoeva@kfs.oeaw.ac.at


{\small
\begin{thebibliography}{150}

\bibitem{CCS} P.G. Casazza, O. Christensen, D.T. Stoeva, Frame expansions in separable Banach spaces, J. Math.
Anal. Appl. 307 (2005) 710--723.

\bibitem{DS} N. Dunford, J. T. Schwartz,
Linear Operators, Part I: General Theory, Vol. VII of Pure and Applied Mathematics, Interscience Publishers, Inc., New York,  1958.

\bibitem{KA} L.V. Kantorovich, G.P. Akilov,
Functional Analysis in Normed Spaces, Pergamon press, New York, 1964.

\bibitem{LS} L.A. Lusternik, V.J. Sobolev,
Elements of Functional Analysis, Hindustan Publ. Corporation, Delhi, and
John Wiley $\&$ Sons, 1974.

\bibitem{MV} R. Meise, D. Vogt, Introduction to Functional Analysis,  Clarendon
Press, Oxford, 1997.

\bibitem{ps2} S. Pilipovi\'c, D.T. Stoeva, Series expansions in Fr\'echet spaces and their
duals, construction of Fr\'echet frames, Journal of Approximation Theory 163 (2011) 1729--1747.

\bibitem{pst} S. Pilipovi\'c, D.T. Stoeva, N. Teofanov, Frames for Fr\'echet spaces,
Bull. Cl. Sci. Math. Nat. Sci. Math.  32 (2007) 69--84.

\end{thebibliography}
}
\end{document}